\documentclass[12pt,thmsa]{article}
\usepackage{amsfonts}

\newtheorem{theorem}{Theorem}
\newtheorem{corollary}{Corollary}
\newtheorem{proposition}{Proposition}
\newtheorem{lemma}{Lemma}

\newcommand{\p}{\Bbb{P}}

\newcommand{\e}{\Bbb{E}}

\newcommand{\ind}{\mbox{\rm 1\hspace{-0.04in}I}}

\newcommand{\R}{\mbox{\rm I\hspace{-0.02in}R}}

\newcommand{\ed}{\stackrel{(d)}{=}}
\newcommand{\eqdef}{\stackrel{\mbox{\tiny$($def$)$}}{=}}

\def\QED{\hfill\vrule height 1.5ex width 1.4ex depth -.1ex \vskip20pt}

\begin{document}

\title{The lower envelope of positive self-similar\\
Markov processes.}
\date{\today}
\maketitle
\begin{center}
{\large L. Chaumont\footnote{E-mail: chaumont@ccr.jussieu.fr}  and
J.C. Pardo\footnote{E-mail: pardomil@ccr.jussieu.fr -- Research
supported by a grant from CONACYT
(Mexico).}\vspace*{0.2in}}\\
Laboratoire de Probabilit\'es et Mod\`eles Al\'eatoires,\\
Universit\'e Pierre et Marie Curie\\ 4, Place Jussieu -- 75252 {\sc
Paris Cedex 05.}
\end{center}

\noindent {\it Abstract}: {\footnotesize We establish integral tests
and laws of the iterated logarithm for the lower envelope of
positive self-similar Markov processes at 0 and $+\infty$. Our
proofs are based on the Lamperti representation and time reversal
arguments. These results extend laws of the iterated logarithm
 for Bessel processes due to Dvoretsky and Erd\"os \cite{de},
 Motoo \cite{mo} and Rivero \cite{ri}.}\\


\noindent
{\it A.M.S. Classification}: {\footnotesize 60 G 18, 60 G 51, 60 B 10.}\\

\section{Introduction}
A real self-similar Markov process $X^{(x)}$, starting from $x$ is a
c\`adl\`ag Markov process which fulfills a scaling property, i.e.,
there exists a constant $\alpha>0$ such that for any $k>0$,
\begin{equation}\label{scale}
\mbox{$\left(kX^{(x)}_{k^{-\alpha}t},\,t\ge0\right)\ed
(X_{t}^{(kx)},\,t\ge0).$}
\end{equation}
Self-similar processes often arise in various parts of probability
theory as limit of re-scaled processes. Their properties have been
studied by the early sixties under the impulse of Lamperti's work
\cite{la1}. The Markov property added to self-similarity provides
some interesting features as noted by Lamperti himself in \cite{la}
where the particular case of positive self-similar Markov processes
is studied. These processes are involved for instance in branching
processes and fragmentation theory. In this paper, we will consider
{\it positive} self-similar Markov processes and refer to them as
{\it pssMp}. Some particularly well known examples are transient
Bessel processes, stable subordinators or more generally, stable
L\'evy processes conditioned to stay positive.

The aim of this work is to describe the lower envelope at 0 and at
$+\infty$ of a large class of pssMp throughout integral tests and
laws of the iterated logarithm (LIL for short). A crucial point in
our arguments is the famous Lamperti representation of self-similar
$\R_+$--valued Markov processes. This transformation enables us to
construct the paths of any such process starting from $x>0$, say
$X^{(x)}$, from those of a L\'evy process. More precisely, Lamperti
\cite{la} found the representation
\begin{equation}\label{lamp}
X_t^{(x)}=x\exp\xi_{\tau(tx^{-\alpha})},\;\;\;0\le t\le
x^{-\alpha}I(\xi)\,,
\end{equation}
under $\p_x$, for $x>0$, where
\[\tau_t=\inf\{s:I_s(\xi)\ge t\}\,,\;\;\;I_s(\xi)=
\int_0^s\exp\alpha\xi_u\,du
\,,\;\;\;I(\xi)=\lim_{t\rightarrow+\infty}I_t(\xi)\,,\] and where
$\xi$ is a real L\'evy process which is possibly killed at
independent exponential time. Note that for $t<I(\xi)$, we have the
equality $\tau_t=\int_0^t\big(X_s^{(x)}\big)^{-\alpha}ds$, so that
(\ref{lamp}) is invertible and yields a one to one relation between
the class of pssMp and the one of L\'evy processes.

In this work, we consider pssMp's which drift towards $+\infty$,
i.e. $\lim_{t\rightarrow+\infty} X_t^{(x)}=+\infty$, a.s. and which
fulfills the Feller property on $[0,\infty)$, so that we may define
the law of a pssMp, which we will call $X^{(0)}$, starting from 0
and with the same transition function as $X^{(x)}$, $x>0$. Bertoin
and Caballero \cite{bc} and Bertoin and Yor \cite{by1} proved that
the family of processes $X^{(x)}$ converges, as $x\downarrow0$, in
the sense of finite dimensional distributions towards $X^{(0)}$ if
and only if the underlying L\'evy process $\xi$ in the Lamperti's
representation is such that
\[(H)\;\;\;\;\;\;\;\;\mbox{$\xi$ is non lattice and}\;\;\;
0<m\eqdef\e(\xi_1)\le\e(|\xi_1|)<+\infty\,.\] As proved by Caballero
and Chaumont in \cite{cc}, the latter condition is also a NASC for
the weak convergence of the family $(X^{(x)})$, $x\ge0$ on the
Skohorod's space of c\`adl\`ag trajectories. In the same article,
the authors also provided a path construction of the process
$X^{(0)}$. The entrance law of $X^{(0)}$ has been described in
\cite{bc} and \cite{by1} as follows: for every $t>0$ and for every
measurable function $f:\R_+\rightarrow\R_+$,
\begin{equation}\label{entlaw}
\e\left(f\left(X_t^{(0)}\right)\right)=
\frac1{m}\e\left(I(-\xi)^{-1}f(tI(-\xi)^{-1})\right)\,.
\end{equation}

Several partial results on the lower envelope of $X^{(0)}$ have
already been established before, the oldest of which being due to
Dvoretsky and Erd\"os \cite{de} and Motoo \cite{mo} who studied the
special case of Bessel processes. More precisely, when $X^{(0)}$ is
a Bessel process with dimension $\delta>2$, we have the following
integral test at 0: if $f$ is an increasing function then
\[\p(X^{(0)}_t<f(t),\mbox{i.o., as
$t\rightarrow0$})=\left\{\begin{array}{l}0\\1\end{array}\right.\;\;\;
\mbox{according as}\;\;\;\int_{0+}
\left(\frac{f(t)}t\right)^{\frac{\delta-2}4}\,\frac{dt}t\;
\left\{\begin{array}{l}<\infty\\=\infty\end{array}\right.\,.\] The
time inversion property of Bessel processes, induces the same
integral test for the behaviour at $+\infty$ of $X^{(x)}$, $x\ge0$.
The test for Bessel processes is extended in the section 4 of this
paper to pssMp's such that the upper tail of the law of the
exponential functional $I(-\xi)$ is regularly varying. Our integral
test is then written in terms of the law of this exponential
functional as shown in Theorem \ref{regul}.

Without giving here an exhaustive list of the results which have
been obtained in that direction, we may also cite Lamperti's own
work \cite{la} who used his representation to describe the
asymptotic behaviour of a pssMp starting from $x>0$ in terms of the
underlying L\'evy process. Some cases where the transition function
of $X^{(x)}$ admits some special bounds have also been studied by
Xiao \cite{xi}.

The most recent result concerns increasing pssMp and is due to
Rivero \cite{ri} who proved the following LIL: suppose that $\xi$ is
a subordinator whose Laplace exponent $\phi$ is regularly varying at
infinity with index $\beta\in(0,1)$ and define the function
$\varphi(t)=\phi(\log|\log t|)/\log|\log t|$, $t>e$, then
\[\liminf_{t\downarrow0}\frac{X_t^{(0)}}{(t\varphi(t))^{1/\alpha}}=
\alpha^{\beta/\alpha}(1-\beta)^{(1-\beta)/\alpha}\;\;\mbox{and}
\;\;\liminf_{t\uparrow+\infty}\frac{X_t^{(0)}}{(t\varphi(t))^{1/\alpha}}=
\alpha^{\beta/\alpha}(1-\beta)^{(1-\beta)/\alpha},\;a.s.\] In
Section 5 of this paper, we extend Rivero's result to pssMp's such
that the logarithm of the upper tail of the exponential functional
$I(-\xi)$ is regularly varying at $+\infty$. In Theorem
\ref{logregul}, we give a LIL for the process $X^{(0)}$ at 0 and for
the processes $X^{(x)}$, $x\ge0$ at $+\infty$. Then the lower
envelope has an explicit form in terms of the tail of the law of
$I(-\xi)$.

All the asymptotic results presented in Sections 4 and 5 are
consequences of general integral tests which are stated and proved
in Section 3 and which may actually be applied in other situations
than our 'regular' and 'logregular' cases. If $F_q$ denotes the tail
of the law of the truncated exponential functionals
$\int_0^{\hat{T}_{-q}}\exp-\xi_s\,ds$,
$\hat{T}_{-q}=\inf\{t:\xi_s\le -q\}$, $q\ge0$, then we will show
that in any case, the knowledge of asymptotic behaviour of $F_q$
suffices to describe the lower envelope of the processes $X^{(x)}$,
$x\ge0$.

Section 2 is devoted to preliminary results. We give a path
decomposition of the process $X^{(0)}$ up to its last passage time
under a fixed level. This process, once reversed, corresponds to a
pssMp whose associated L\'evy process in the Lamperti transformation
is $-\xi$. In particular, this allows us to get an expression of the
last passage time process of $X^{(0)}$ in terms of $I(-\xi)$. The
description of the last passage process is then used in section 3
for the study of the lower envelope of $X^{(x)}$.

\section{Time reversal and last passage time of $X^{(0)}$}

We consider processes defined on the space ${\cal D}$ of c\`adl\`ag
trajectories on $[0,\infty)$, with real values. The space ${\cal D}$
is endowed with the Skorohod's topology and $\p$
will be our reference probability measure.\\

{\it In all the rest of the paper, $\xi$ will be a L\'evy process
satisfying condition $(H)$. With no loss of generality, we will also
suppose that $\alpha=1$. Indeed, we see from (\ref{scale}) that if
$X^{(x)}$, $x\ge0$, is a pssMp with index $\alpha>0$, then
$\big(X^{(x)}\big)^{\alpha}$ is a pssMp with index 1. Therefore, the
integral tests and LIL established in the sequel can easily be
interpreted for any $\alpha>0$.}\\

Let us define the family of positive self-similar Markov processes
$\hat{X}^{(x)}$ whose Lamperti's representation is given by
\begin{equation}\label{lamp2}
\hat{X}^{(x)}=\left(x\exp\hat{\xi}_{\hat{\tau}(t/x)},\, 0\le t\le
xI(\hat{\xi})\right)\,,\;\;\;x>0\,,\end{equation} where
$\hat{\xi}=-\xi$,
$\hat{\tau}_t=\inf\{s:\int_0^s\exp\hat{\xi}_u\,du\ge t\}$, and
$I(\hat{\xi})=\int_0^\infty\exp\hat{\xi_s}\,ds$. We emphasize that
the r.v. $xI(\hat{\xi})$, corresponds to the first time at which the
process $\hat{X}^{(x)}$ hits 0, i.e.
\begin{equation}\label{expfunc}
xI(\hat{\xi})=\inf\{t:\hat{X}^{(x)}_t=0\}\,,
\end{equation}
moreover, for each $x>0$, the process $\hat{X}^{(x)}$ hits 0
continuously, i.e. $\hat{X}^{(x)}(xI(\hat{\xi})-)=0$.

We now fix a decreasing sequence $(x_n)$, $n\ge1$ of positive real
numbers which tends to $0$ and we set
\[U(y)=\sup\{t:X^{(0)}_t\le y\}\,.\]
The aim of this section is to establish a path decomposition of the
process $X^{(0)}$ reversed at time $U(x_1)$ in order to get a
representation of this time in terms of the exponential functional
$I(\hat{\xi})$, see Corollaries \ref{cor1} and \ref{cor3} below.

To simplify the notations, we set $\Gamma=X_{U(x_1)-}^{(0)}$ and we
will denote by $K$ the support of the law of $\Gamma$. We will see
in Lemma \ref{overshoot} that actually $K=[0,x_1]$. For any process
$X$ that we consider here, we make the convention that $X_{0-}=X_0$.

\begin{proposition}\label{ret} The  law of the process
$\hat{X}^{(x)}$ is a regular version of the law of
 the process
$$\hat{X}\eqdef(X^{(0)}_{(U(x_1)-t)-},\,0\le t\le U(x_1))\,,$$
conditionally on $\Gamma=x$, $x\in K$.
\end{proposition}
{\it Proof}: The result is a consequence of Nagasawa's theory of
time reversal for Markov processes. First, it follows from Lemma 2
in \cite{by1} that the resolvent operators of $X^{(x)}$ and
$\hat{X}^{(x)}$, $x>0$ are in duality with respect to the Lebesgue
measure. More specifically, for every $q\ge0$, and measurable
functions $f,g:(0,\infty)\rightarrow\R_+$, with
\[V^qf(x)\eqdef\e\left(\int_0^\infty e^{-qt}f(X_t^{(x)})\,dt\right),\;\;\;
\mbox{and}\;\;\;\;\hat{V}^qf(x)\eqdef\e\left(\int_0^\zeta
e^{-qt}f(\hat{X}_t^{(x)})\,dt\right),\] we have
\begin{equation}\label{scal1}\int_0^\infty f(x)\hat{V}^qg(x)\,dx=\int_0^\infty
g(x)V^qf(x)\,dx\,.\end{equation} Let $p_t(dx)$ be the entrance law
of $X^{(0)}$ at time $t$, then it follows from the scaling property
that for any $t>0$, $p_t(dx)=p_{1}(dx/t)$, hence $ \int_0^\infty
p_t(dx)\,dt=\int_0^\infty p_1(dy)/y\,dx$ for all $x>0$, where from
(\ref{entlaw}), $\int_0^\infty p_1(dy)/y\,dy=m^{-1}$. In other
words, the resolvent measure of $\delta_{\{0\}}$ is proportional to
the Lebesgue measure, i.e.:
\begin{equation}\label{scal2}m^{-1}\int_0^\infty
f(x)\,dx=\e\left(\int_0^\infty f(X_t^{(0)})\,dt\right)\,.
\end{equation}
Conditions of Nagasawa's theorem are satisfied as shown in
(\ref{scal1}) and (\ref{scal2}), then it remains to apply this
result to $U(x_1)$ which is a return time such that
$\p(0<U(x_1)<\infty)=1$, and the proposition is proved.\QED

\noindent Another way to state Proposition \ref{ret} is to say that
for any $x\in K$, the returned process
$(\hat{X}_{(xI(\hat{\xi})-t)-},\,0\le t\le xI(\hat{\xi}))$, has the
same law as $(X_t^{(0)},0\le t<U(x_1))$ given $\Gamma=x$. In
\cite{by1}, the authors show that when the semigroup operator of
$X^{(0)}$ is absolutely continuous with respect to the Lebesgue
measure with density $p_t(x,y)$, this process is an $h$-process of
$X^{(0)}$, the corresponding harmonic function being
$h(x)=\int_0^\infty p_t(x,1)\,dt$.

For $y>0$, we set
\[\hat{S}_y=\inf\{t:\hat{X}_t\le y\}\,.\]

\begin{corollary}\label{cor2}
Between the passage times $\hat{S}_{x_n}$ and $\hat{S}_{x_{n+1}}$,
the process $\hat{X}$ may be described as follows:
\[\left(\hat{X}_{\hat{S}(x_n)+t},\,0\le t\le
\hat{S}_{x_{n+1}}-\hat{S}_{x_n}\right)=\left(\Gamma_n
\exp\hat{\xi}^{(n)}_{\hat{\tau}^{(n)}(t/\Gamma_n)},\, 0\le t\le
H_n\right),\;\;n\ge1,\] where the processes $\hat{\xi}^{(n)}$,
$n\ge1$ are independent between themselves and have the same law as
$\hat{\xi}$. Moreover the sequence $(\hat{\xi}^{(n)})$ is
independent of $\Gamma$ defined above and
\[\left\{\begin{array}{l}
\hat{\tau}_t^{(n)}=\inf\{s:\int_0^s\exp\hat{\xi}_u^{(n)}\,du\ge t\}\\
H_n=\Gamma_n\int_0^{\hat{T}^{(n)}(\log(x_{n+1}/\Gamma_n))}\exp\hat{\xi}^{(n)}_s\,ds\\
\Gamma_{n+1}=\Gamma_{n}\exp\hat{\xi}^{(n)}_{\hat{T}^{(n)}(\log
x_{n+1}/\Gamma_{n})},\,n\ge1,\;\;\;\Gamma_1=\Gamma\\
\hat{T}_z^{(n)}=\inf\{t:\hat{\xi}^{(n)}_t\le z\}.
\end{array}\right.\]
For each $n$, $\Gamma_n$ is independent of $\xi^{(n)}$ and
\begin{equation}\label{mmloi}
x_n^{-1}\Gamma_n\ed x_1^{-1}\Gamma\,.
\end{equation}
\end{corollary}
{\it Proof}: From (\ref{lamp2}) and Proposition \ref{ret}, the
process $\hat{X}$ may be described as
\[\hat{X}=\left(\Gamma\exp\hat{\xi}^{(1)}_{\hat{\tau}^{(1)}(t/\Gamma)},\,0\le t\le
U(x_1)\right)\,,\] where $\hat{\xi}^{(1)}\ed\hat{\xi}$ is
independent of $\Gamma=X^{(0)}_{U(x_1)-}$ and
$\hat{\tau}_t^{(1)}=\inf\{s:\int_0^s\exp\hat{\xi}_u^{(1)}\,du\ge
t\}$. Note that $\Gamma\le x_1$, a.s., so between the passages times
$\hat{S}_{x_1}=0$ and $\hat{S}_{x_2}$, the process $\hat{X}$ is
clearly described as in the statement with
$\hat{\xi}^{(1)}=\hat{\xi}$ and $\hat{S}_{x_2}-\hat{S}_{x_1}=H_1=
\Gamma\int_0^{\hat{T}^{(1)}(\log(x_{2}/\Gamma))}\exp\hat{\xi}_s^{(1)}\,ds$.

Now if we set
$\hat{\xi}^{(2)}\eqdef(\hat{\xi}^{(1)}_{\hat{T}^{(1)}(\log
x_2/\Gamma_1)+t}-\hat{\xi}^{(1)}_{\hat{T}^{(1)}(\log
x_2/\Gamma_1)},\,t\ge0)$, then with the definitions of the
statement,
\begin{eqnarray}
&&(\hat{X}_{\hat{S}(x_2)+t},\,t\ge0)=(\Gamma_2
\exp\hat{\xi}^{(2)}_{\hat{\tau}^{(2)}(t/\Gamma_2)},\,t\ge0)\;\;\;\;\mbox{and}\label{ord2}\\
&&\hat{S}_{x_{3}}-\hat{S}_{x_2}=\inf\{t:\hat{X}_{\hat{S}(x_2)+t}\le
x_3\}=H_2\,.\nonumber
\end{eqnarray}
The process $\hat{\xi}^{(2)}$ is independent of
$[(\hat{\xi}^{(1)}_t,\,0\le t\le \hat{T}^{(1)}(\log
x_2/\Gamma_1)),\,\Gamma_1]$, hence it is clear that we do not change
the law of $\hat{X}$ if, by reconstructing it according to this
decomposition, we replace $\hat{\xi}^{(2)}$ by a process with the
same law which is independent of $[\hat{\xi}^{(1)},\Gamma_1]$.
Moreover, $\hat{\xi}^{(2)}$ is independent of $\Gamma_2$. Relation
(\ref{mmloi}) is a consequence of the scaling property. Indeed, we
have
\[\left(\frac{x_2}{x_1}X^{(0)}_{tx_1/x_2},\,0\le
t\le \frac{x_2}{x_1}U(x_1)\right)\ed\left(X^{(0)}_{t},\,0\le t\le
U(x_2)\right)\,,\] which implies the identities in law
\begin{equation}\label{id1}
x_1^{-1}X^{(0)}_{U(x_1)-}\ed
x_2^{-1}X^{(0)}_{U(x_2)-}\,,\;\;\;\mbox{and}
\,\;\;\;x_1^{-1}U(x_1)\ed x_2^{-1}U(x_2)\,.
\end{equation} On the other hand, we see from the definition of $\hat{X}$ in
Proposition \ref{ret} that
\[\left(\hat{X}_{\hat{S}(x_2)+t},\,0\le t\le U(x_1)-\hat{S}(x_2)\right)=
\left(X^{(0)}_{(U(x_2)-t)-},\,0\le t\le U(x_2)\right)\,.\] Then, we
obtain (\ref{mmloi}) for $n=2$ from this identity, (\ref{ord2}) and
(\ref{id1}). The proof follows by induction.\QED

\begin{corollary}\label{cor1}  With the same notations as in
Corollary $\ref{cor2}$, the time $U(x_n)$ may be decomposed into the
sum
\begin{equation}\label{last2}
U(x_n)=\sum_{k\ge
n}\Gamma_k\int_0^{\hat{T}^{(k)}(\log(x_{k+1}/\Gamma_k))}
\exp\hat{\xi}^{(k)}_s\,ds\,,\;\;\;a.s.
\end{equation} In particular, for all $z_n>0$, we have
\begin{equation}\label{last3}
z_{n}\ind_{\{\Gamma_n\ge
z_n\}}\int_0^{\hat{T}^{(n)}(\log(x_{n+1}/z_n))}
\exp\hat{\xi}^{(n)}_s\,ds \le U(x_n)\le x_n
I(\overline{\xi}^{(n)})\,,\;\;\;a.s.,
\end{equation}
where $\overline{\xi}^{(n)}$, $n\ge1$ are L\'evy processes with the
same law as $\hat{\xi}$.
\end{corollary}
{\it Proof}: Identity (\ref{last2}) is a consequence of Corollary
\ref{cor2} and the fact that $U(x_n)=\sum_{k\ge
n}\hat{S}_{k+1}-\hat{S}_k$. The first inequality in (\ref{last3}) is
a consequence of (\ref{last2}), which implies:
$\Gamma_n\int_0^{\hat{T}^{(n)}(\log(x_{n+1}/\Gamma_n))}
\exp\hat{\xi}^{(n)}_s\,ds\le U(x_n)$.

To prove the second inequality in (\ref{last3}), it suffices to note
that by Proposition \ref{ret} and the strong Markov property at time
$\hat{S}(x_n)$, for any $n\ge1$, we have the representation
\[\left(\hat{X}_{\hat{S}(x_n)+t},\,0\le t\le U(x_1)-\hat{S}(x_n)\right)=
\left(\Gamma_n\exp\overline{\xi}^{(n)}_{\overline{\tau}^{(n)}(t/\Gamma_n)},\,
0\le t\le U(x_1)-\hat{S}(x_n)\right)\,,\] where
$\overline{\tau}^{(n)}_t=\inf\{s:\int_0^s\exp\overline{\xi}^{(n)}_u\,du>t\}$
and $\Gamma_n=\hat{X}_{\hat{S}(x_n)}$ (see in Corollary \ref{cor2})
is independent of $\overline{\xi}^{(n)}$ which has the same law as
$\hat{\xi}$. It remains to note from (\ref{expfunc}) that
$U(x_1)-\hat{S}(x_n)=U(x_n)=\Gamma_nI(\overline{\xi}^{(n)})$ and
that $\Gamma_n\le x_n$.\QED

\noindent To establish our asymptotic results at $+\infty$, we will
also need to estimate the law of the time $U(x)$ when $x$ is large.
The same reasoning as we did for a sequence which tends to 0 can be
done for a sequence which tends to $+\infty$ as we show in the
following result.

\begin{corollary}\label{cor3}  Let $(y_n)$ be an increasing sequence of positive real
numbers which tends to $+\infty$. There exists some sequences
$(\check{\xi}^{(n)})$, $(\tilde{\xi}^{(n)})$ and
$(\check{\Gamma}_n)$, such that for each $n$,
$\check{\xi}^{(n)}\ed\tilde{\xi}^{(n)}\ed\hat{\xi}$,
$\check{\Gamma}_n\ed\Gamma$, $\check{\Gamma}_n$ and
$\check{\xi}^{(n)}$ are independent; moreover the L\'evy processes
$(\check{\xi}^{(n)})$ are independent between themselves and we have
for all $z_n>0$,
\begin{equation}\label{last4}
z_{n}\ind_{\{\check{\Gamma}_n\ge
z_n\}}\int_0^{\check{T}^{(n)}(\log(y_{n-1}/z_n))}
\exp\check{\xi}^{(n)}_s\,ds \le U(y_n)\le y_n
I(\tilde{\xi}^{(n)})\,,\;\;\;a.s.
\end{equation}
where $\check{T}^{(n)}_z=\inf\{t:\check{\xi}^{(n)}_t\le z\}$.
\end{corollary}
{\it Proof}: Fix an integer $n\ge1$ and define the decreasing
sequence $x_1,\dots,x_n$ by $x_n=y_1, x_{n-1}=y_2,\dots,x_1=y_n$,
then  construct the sequences
$\hat{\xi}^{(1)},\dots,\hat{\xi}^{(n)}$ and
$\Gamma_1,\dots,\Gamma_n$ from $x_1,\dots,x_n$ as in Corollary
\ref{cor2} and construct the sequence
$\overline{\xi}^{(1)},\dots,\overline{\xi}^{(n)}$ as in Corollary
\ref{cor1}. Now define
$\check{\xi}^{(1)}=\hat{\xi}^{(n)},\check{\xi}^{(2)}=
\hat{\xi}^{(n-1)},\dots,\check{\xi}^{(n)}=\hat{\xi}^{(1)}$ and
$\tilde{\xi}^{(1)}=\overline{\xi}^{(n)},\tilde{\xi}^{(2)}=
\overline{\xi}^{(n-1)},\dots,\tilde{\xi}^{(n)}=\overline{\xi}^{(1)}$
and $\check{\Gamma}_{1}=\Gamma_n,\check{\Gamma}_{2}=
\Gamma_{n-1},\dots,\check{\Gamma}_{n}=\Gamma_{1}$. Then from
(\ref{last3}), we deduce that for any $k=2,\dots,n$,
\[
z_{k}\ind_{\{\check{\Gamma}_k\ge
z_k\}}\int_0^{\check{T}^{(k)}(\log(y_{k-1}/z_k))}
\exp\check{\xi}^{(k)}_s\,ds \le U(y_k)\le y_k
I(\tilde{\xi}^{(k)})\,,\;\;\;a.s.
\]
Hence the whole sequences $(\tilde{\xi}^{(n)})$,
$(\check{\xi}^{(n)})$ and $(\check{\Gamma}_n)$ are well constructed
and fulfill the desired properties. \QED

\noindent {\bf Remark}: We emphasize that
$\hat{T}^{(n)}(\log(x_{n+1}/\Gamma_n))=0$, a.s. on the event
$\Gamma_n\le x_{n+1}$; moreover, we have $\Gamma_n\le x_n$, a.s., so
the first inequality in $(\ref{last3})$ is relevant only when
$x_{n+1}<z_n<x_n$. Similarly, in Corollary \ref{cor1}, the first
inequality in $(\ref{last4})$ is relevant only when
$y_{n-1}<z_n<y_n$.\\

\noindent We end this section with the computation of the law of
$\Gamma$. Recall that the upward ladder height process $(\sigma_t)$
associated to $\xi$ is the subordinator which corresponds to the
right continuous inverse of the local time at 0 of the reflected
process $(\xi_t-\sup_{s\le t}\xi_s)$, see \cite{be} Chap. V for a
proper definition. We denote by $\nu$ the L\'evy measure of
$(\sigma_t)$.

\begin{lemma}\label{overshoot} The law of $\Gamma$ is characterized
as follows:
\[\log x_1^{-1}\Gamma\ed -{\cal U}Z\,,\]
where ${\cal U}$ and $Z$ are independent r.v.'s, $U$ is uniformly
distributed over $[0,1]$ and the law of $Z$ is given by:
\begin{equation}\label{over}
\p(Z>u)=\e(\sigma_1)^{-1}\int_{(u,\infty)}s\,\nu(ds),\;\;\;u\ge0\,.
\end{equation}
In particular, for all $\eta<x_1$, $\p(\Gamma>\eta)>0$.
\end{lemma}
{\it Proof}. It is proved in \cite{dm} that under our hypothesis,
(that is $\e(|\hat{\xi}_1|)<+\infty$, $E(\hat{\xi}_1)<0$ and $\xi$
is not arithmetic), the overshoot process of $\xi$ converges in law,
that is
\[\hat{\xi}_{\hat{T}(x)}-x\longrightarrow
-{\cal U}Z,\;\;\;\mbox{in law as $x$ tends to $-\infty$,}\] and the
limit law is computed in \cite{cho} in terms of the upward ladder
height process $(\sigma_t)$.

On the other hand, we proved in Corollary \ref{cor2}, that
\begin{eqnarray*}
x_{n+1}^{-1}\Gamma_{n+1}&=&\exp[\hat{\xi}^{(n)}_{\hat{T}^{(n)}(\log
x_{n+1}/\Gamma_{n})}-\log x_{n+1}/\Gamma_n]\ed x_1^{-1}\Gamma\\
&\ed&\exp[\hat{\xi}_{\hat{T}(\log x_{n+1}/x_n+\log
x_1^{-1}\Gamma)}-\log x_{n+1}/x_n-\log x_1^{-1}\Gamma]\,.
\end{eqnarray*}
Then by taking $x_n=e^{-n^2}$, we deduce from these equalities that
$\log x_1^{-1}\Gamma$ has the same law as the limit overshoot of the
process $\hat{\xi}$, i.e.
\[\hat{\xi}_{\hat{T}(x)}-x\longrightarrow
\log x_1^{-1}\Gamma,\;\;\;\mbox{in law as $x$ tends to $-\infty$.}\]
\QED

\noindent As a consequence of the above results we have the
following identity in law:
\[U(x)\ed \frac xx_1\Gamma I(\hat{\xi})\,,\]
($\Gamma$ and $I(\hat{\xi})$ being independent) which has been
proved in \cite{bc}, Proposition 3 in the special case where the
process $X^{(0)}$ is increasing.

\section{The lower envelope}
\setcounter{equation}{0}

The main result of this section consists in integral tests at 0 and
$+\infty$ for the lower envelope of $X^{(0)}$. When no confusion is
possible, we set $I\eqdef
I(\hat{\xi})=\int_0^{\infty}\exp\hat{\xi}_s\,ds$. This theorem means
in particular that the asymptotic behaviour of $X^{(0)}$ only
depends on the tail behaviour of the law of $I$, and on this of the
law of $\int_0^{\hat{T}_{-q}}\exp\hat{\xi}_s\,ds$, with
$\hat{T}_{x}=\inf\{t:\hat{\xi}_t\le x\}$, for $x\le 0$. So also we
set
\[I_q\eqdef\int_0^{\hat{T}_{-q}}\exp\hat{\xi}_s\,ds\,,\;\;\;\;
F(t)\eqdef\p(I>t)\,,\;\;\;\;F_q(t)\eqdef\p(I_q>t)\,.\] The following
lemma will be used to show that actually, in many particular cases,
$F$ suffices to describe the envelope of $X^{(0)}$.
\begin{lemma}\label{lem3}
Assume that there exists $\gamma>1$ such that,
$\limsup_{t\rightarrow+\infty}F(\gamma t)/F(t)<1$. For any $q>0$ and
$\delta>\gamma e^{-q}$,
\[\liminf_{t\rightarrow+\infty}\frac{F_q((1-\delta)t)}{F(t)}>0\,.\]
\end{lemma}
{\it Proof}: It follows from the decomposition of $\xi$ into the two
independent processes $(\hat{\xi}_s,\,s\le \hat{T}_{-q})$ and
$\hat{\xi}'\eqdef(\hat{\xi}_{s+\hat{T}_{-q}}-
\hat{\xi}_{\hat{T}_{-q}},\,s\ge0)$ that
\[I=I_q+e^{\hat{\xi}_{\hat{T}_{-q}}}I'\le I_q+e^{-q}I'\]
where $I'=\int_0^\infty\exp\hat{\xi}'_s\,ds$ is a copy of $I$ which
is independent of $I_q$. Then we can write for any $q>0$ and
$\delta\in(0,1)$, the inequalities
\begin{eqnarray*}
\p(I>t)&\le&\p(I_q+e^{-q}I'\ge t)\\
&\le&\p(I_q>(1-\delta)t)+\p(e^{-q}I>\delta t)\,,
\end{eqnarray*}
so that if moreover, $\delta>\gamma e^{-q}$ then
\[1-\frac{\p(I>\gamma t)}{\p(I>t)}\le1-
\frac{\p(I>e^q\delta t)}{\p(I>t)}\le
\frac{\p(I_q>(1-\delta)t)}{\p(I>t)}\,.\] \QED

\noindent We start by stating the integral test at time 0.

\begin{theorem}\label{main1}
The lower envelope of $X^{(0)}$
at $0$ is described as follows:\\

\noindent Let $f$ be an increasing function.
\begin{itemize}
\item[$(i)$] If
\[\int_{0+}F\left(\frac t{f(t)}\right)\,\frac{dt}t<\infty\,,\]
then for all $\varepsilon>0$,
\[\p(X^{(0)}_t<(1-\varepsilon)f(t),\mbox{i.o., as
$t\rightarrow0$})=0\,.\] \item[$(ii)$] If for all $q>0$,
\[\int_{0+}F_q\left(\frac t{f(t)}\right)\,\frac{dt}t=\infty\,,\]
then for all $\varepsilon>0$,
\[\p(X^{(0)}_t<(1+\varepsilon)f(t),\mbox{i.o., as
$t\rightarrow0$})=1\,.\] \item[$(iii)$] Suppose that $t\mapsto
f(t)/t$ is increasing. If there exists $\gamma>1$ such that,
\[\mbox{$\limsup_{t\rightarrow+\infty}\p(I>\gamma t)/\p(I>t)<1$ and if}\;\;\;
\int_{0+}F\left(\frac t{f(t)}\right)\,\frac{dt}t=\infty\,,\] then
for all $\varepsilon>0$,
\[\p(X^{(0)}_t<(1+\varepsilon)f(t),\mbox{i.o., as
$t\rightarrow0$})=1\,.\]
\end{itemize}
\end{theorem}

\noindent{\it Proof}: Let $(x_n)$ be a decreasing sequence such that
$\lim_nx_n=0$. Recall the notations of Section 2. We define the
events \[A_n=\{\mbox{There exists $t\in[U(x_{n+1}),U(x_n)]$ such
that $X_t^{(0)}<f(t)$.}\}\,.\] Since $U(x_n)$ tends to 0, a.s. when
$n$ goes to $+\infty$,  we have:
\begin{equation}\label{limsup}\{X_t^{(0)}<f(t),\mbox{i.o., as
$t\rightarrow0$}\}=\limsup_nA_n\,.
\end{equation} Since $f$ is increasing, the following inclusions hold:
\begin{equation}\label{inc1}\{x_n\le
f(U(x_n))\}\subset A_n\subset \{x_{n+1}\le
f(U(x_n))\}\,.\end{equation}

\noindent Then we prove the convergent part $(i)$. Let us choose
$x_n=r^{-n}$ for $r>1$, and recall from relation (\ref{last3}) above
that $U(r^{-n})\le r^{-n} I(\overline{\xi}^{(n)})$. From this
inequality and (\ref{inc1}), we can write:
\begin{equation}\label{inc5}
A_n\subset \{r^{-(n+1)}\le
f(r^{-n}I(\overline{\xi}^{(n)}))\}\,.\end{equation} Let us denote
$I(\hat{\xi})$ simply by $I$. From Borel Cantelli's Lemma,
(\ref{inc5}) and (\ref{limsup}),
\begin{equation}\label{bc1} \mbox{if
$\sum_n\p(r^{-(n+1)}\le f(r^{-n}I))<\infty$ then $\p(X_t^{(0)}<f(t),
\mbox{i.o., as $t\rightarrow0$})=0$.}
\end{equation}
Note that $\int_{1}^{+\infty}\p(r^{-t}\le
f(r^{-t}I))\,dt=\int_{0+}^{+\infty}\p(s<f(s)I,\,s<I/r)/(s\log
r)\,ds$, hence since $f$ is increasing, we have the inequalities:
\begin{equation}\label{ineq}\sum_{n=1}^\infty\p(r^{-n}\le f(r^{-(n+1)}I))\le
\int_{0+}^{+\infty}\p\left(\frac
s{f(s)}<I,\,s<\frac{I}r\right)\,\frac{ds}{s\log r}\le
\sum_{n=1}^\infty\p(r^{-(n+1)}\le f(r^{-n}I))\,.\end{equation} With
no loss of generality, we can restrict ourself to the case $f(0)=0$,
so it is not difficult to check that for any $r>1$,
\begin{equation}\label{ineq2}\int_{0+}\p\left(\frac
s{f(s)}<I,\,s<\frac{I}r\right)\,\frac{ds}s<+\infty,\;\;\;\mbox{if
and only if}\;\;\;\int_{0+}\p\left(\frac
s{f(s)}<I\right)\,\frac{ds}s<+\infty\,.\end{equation} Suppose the
latter condition holds, then from (\ref{ineq}), for all $r>1$,
$\sum_{n=2}^\infty\p(r^{-(n+1)}\le r^{-2}f(r^{-n}I))<+\infty$ and
from (\ref{bc1}), for all $r>1$, $\p(X_t^{(0)}<r^{-2}f(t),
\mbox{i.o., as $t\rightarrow0$})=0$
which proves the desired result.\\

\noindent Now we prove the divergent part $(ii)$. Again, we choose
$x_n=r^{-n}$ for $r>1$, and $z_n=kr^{-n}$, where
$k=1-\varepsilon+\varepsilon/r$ and $0<\varepsilon<1$, (so that
$x_{n+1}<z_n<x_n$). We set
$$B_n=\{r^{-n}\le f_{r,\varepsilon}(kr^{-n}\ind_{\{\Gamma_n\ge
kr^{-n}\}}{I}^{(n)})\}\,,$$ where, $f_{r,\varepsilon}(t)=rf(t/k)$
and with the same notations as in Corollary \ref{cor1}, for each
$n$,
\begin{equation}\label{ids2}I^{(n)}\eqdef
\int_0^{\hat{T}^{(n)}(\log(x_{n+1}/z_n))}
\exp\hat{\xi}^{(n)}_s\,ds\ed\int_0^{\hat{T}(\log(1/rk))}
\exp\hat{\xi}_s\,ds\end{equation}
 is independent of $\Gamma_n$, and $\Gamma_n$ is such that
$x^{-1}_n\Gamma_n\ed x^{-1}_1\Gamma$. Moreover the r.v.'s $I^{(n)}$,
$n\ge 1$ are independent between themselves and identity
(\ref{ids2}) shows that they have the same law as $I_q$ defined in
Lemma \ref{lem3}, where $q=-\log(1/rk)$. With no loss of generality,
we may assume that $f(0)=0$, so that we can write $B_n=\{r^{-n}\le
f_{r,\varepsilon}(kr^{-n}{I}^{(n)}),\,\Gamma_n\ge kr^{-n}\}$ and
from the above arguments we deduce
\begin{equation}\label{eq2}
\p(B_n)=\p(r^{-n}\le f_{r,\varepsilon}(kr^{-n}{I}_q))\p(\Gamma\ge
kr^{-1})\,.
\end{equation}
The arguments which are developed above to show (\ref{ineq}) and
(\ref{ineq2}), are also valid if we replace $I$ by $I_q$. Hence from
the hypothesis, since $\int_{0+}\p(s<f(s)I_q)\,ds/s=+\infty$, then
from (\ref{ineq}) and (\ref{ineq2}) applied to $I_q$, we have
$\sum_{n=1}^\infty\p(r^{-(n+1)}\le
f(r^{-n}I_q))=\sum_{n=1}^\infty\p(r^{-n}\le
f_{r,\varepsilon}(kr^{-n}I_q))=\infty$, and from (\ref{eq2}) we have
$\sum_n\p(B_n)=+\infty$. Then another application of (\ref{eq2}),
gives for any $n$ and $m$,
\begin{eqnarray*}
\p(B_n\cap B_m)&\le&\p(r^{-n}\le
f_{r,\varepsilon}(kr^{-n}{I}_q))\p(r^{-m}\le
f_{r,\varepsilon}(kr^{-m}{I}_q))\\
\p(B_n\cap B_m)&\le&\p(\Gamma\ge kr^{-1})^{-2}\p(B_n)\p(B_m)\,,
\end{eqnarray*}
where $\p(\Gamma\ge kr^{-1})>0$, from (\ref{over}). Hence from the
extension of Borel-Cantelli's lemma given in \cite{ks},
\begin{equation}\label{eq3}\p(\limsup B_n)\ge \p(\Gamma\ge kr^{-1})^{2}>0\,.
\end{equation}
Then recall from Corollary \ref{cor1} the inequality
$kr^{-n}\ind_{\{\Gamma_n\ge kr^{-n}\}}I^{(n)} \le U(r^{-n})$ which
implies from (\ref{inc1}) that $B_n\subset A_n$, (where in the
definition of $A_n$ we replaced $f$ by $f_{r,\varepsilon}$). So,
from (\ref{eq3}), $\p(\limsup_n A_n)>0$, but since $X^{(0)}$ is a
Feller process and since $\limsup_nA_n$ is a tail event, we have
$\p(\limsup_n A_n)=1$. We deduce from the scaling property of
$X^{(0)}$ and (\ref{limsup}) that
\begin{eqnarray*}
\p(X_t^{(0)}\le f_{r,\varepsilon}(t),\mbox{i.o., as
$t\rightarrow0$.})&=&\p(X_{kt}^{(0)}\le
rf(t),\mbox{i.o., as $t\rightarrow0$.})\\
&=&\p(X_t^{(0)}\le k^{-1}rf(t),\mbox{i.o., as
$t\rightarrow0$.})=1\,.
\end{eqnarray*}
Since $k=1-\varepsilon+\varepsilon/r$, with $r>1$ and
$0<\varepsilon<1$ arbitrary chosen, we obtain $(ii)$.\\

\noindent  Now we prove the divergent part $(iii)$. The sequences
$(x_n)$ and $(z_n)$ are defined as in the proof of $(ii)$ above.
Recall that $q=-\log (1/rk)$ and take $\delta>\gamma e^{-q}$ as in
Lemma \ref{lem3}. With no loss of generality, we may assume that
$f(t)/t\rightarrow0$, as $t\rightarrow0$. Then from the hypothesis
in $(iii)$ and Lemma \ref{lem3}, we have
\[\int_{0+}F_q\left(\frac {(1-\delta)t}{f(t)}\right)\,\frac{dt}t=\infty\,.\]
As already noticed above, this is equivalent to
$\int_{1}^{+\infty}\p((1-\delta)r^{-t}\le f(r^{-t}I_q))\,dt=\infty$.
Since $t\mapsto f(t)/t$ increases,
$\int_{1}^{+\infty}\p((1-\delta)r^{-t}\le
f(r^{-t}I_q))\,dt\le\sum_1^\infty \p((1-\delta)r^{-n}\le
f(r^{-n}I_q))=\infty$. Set
$f_r^{(\delta)}(t)=(1-\delta)^{-1}f(t/k)$, then
$$\sum_1^\infty \p(r^{-n}\le f_r^{(\delta)}(kr^{-n}I_q))=\infty\,.$$
Similarly as in the proof of $(ii)$, define $B'_n=\{r^{-n}\le
f_{r}^{(\delta)}(kr^{-n}{I}^{(n)}),\,\Gamma_n\ge kr^{-n}\}$. Then
$B'_n\subset A_n$, (where in the definition of $A_n$ we replaced $f$
by $f_{r}^{(\delta)}$). From the same arguments as above, since
$\sum^\infty\p(B'_n)=\infty$, we have $\p(\limsup_n A_n)=1$, hence
from the scaling property of $X^{(0)}$ and (\ref{limsup})
\begin{eqnarray*}
\p(X_t^{(0)}\le f_{r}^{(\delta)}(t),\mbox{i.o., as
$t\rightarrow0$.})&=&\p(X_{kt}^{(0)}\le
(1-\delta)^{-1}f(t),\mbox{i.o., as $t\rightarrow0$.})\\
&=&\p(X_t^{(0)}\le k^{-1}(1-\delta)^{-1}f(t),\mbox{i.o., as
$t\rightarrow0$.})=1\,.
\end{eqnarray*}
Since $k=1-\varepsilon+\varepsilon/r$, with $r>1$ and
$0<\varepsilon<1$ and $\delta>\gamma
e^{-q}=\gamma/(r+\varepsilon(1-r))$, by choosing $r$ sufficiently
large and $\varepsilon$ sufficiently small, $\delta$ can be taken
sufficiently small so that $k^{-1}(1-\delta)^{-1}$ is arbitrary
close to 1.\QED

\noindent The divergent part of the integral test at $+\infty$
requires the following Lemma.

\begin{lemma}\label{finite}
For any L\'evy process $\xi$ such that
$0<\e(\xi_1)\le\e(|\xi_1|)<\infty$, and for any $q\ge0$,
\[\mbox{$\e\left(\left|\inf_{t\le
{T}_{q}}{\xi}_t\right|\right)<\infty$}\,,\] where
$T_q=\inf\{t:\xi_t\ge q\}$.\end{lemma} {\it Proof}. The proof bears
upon a result on stochastic bounds for L\'evy processes due to Doney
\cite{do} which we briefly recall. Let $\nu_n$ be the time at which
the $n$-th jump of $\xi$ whose value lies in $[-1,1]^c$, occurs and
define
\[I_n=\inf_{\nu_n\le t<\nu_{n+1}}\xi_t\,.\]
Theorem 1.1 in \cite{do} asserts that the sequence $(I_n)$ admits
the representation
\[I_n=S^{(-)}_n+\tilde{\mbox{\em\i}}_0,\,\;\;n\ge0\,,\]
where $S^{(-)}$ is a random walk with the same distribution as
$(\xi(\nu_n),\,n\ge0)$ and $\tilde{\mbox{\em\i}}_0$ is independent
of $S^{(-)}$. For $a\ge0$, let $\sigma(a)=\min\{n:S^{(-)}_n>a\}$,
then for any $q\ge0$, we have the inequality
\begin{equation}\label{1}
\min_{n\le\sigma(q+|\tilde{\mbox{\em\i}}_0|)}(S^{(-)}_n+\tilde{\mbox{\em\i}}_0)
\le \inf_{t\le T_q}\xi_t\,.
\end{equation}
 On the other hand, it
follows from our hypothesis on $\xi$ that $0<\e(S_1^{(-)})\le
\e(|S^{(-)}_1|)<+\infty$, hence from Theorem 2 of \cite{ja} and its
proof, there exists a finite constant $C$ which depends only on the
law of $S^{(-)}$ such that for any $a\ge0$,
\begin{equation}\label{2}
\mbox{$\e\left(\left|\min_{n\le\sigma(a)}S^{(-)}_n\right|\right)\le
C\e(\sigma(a))\e(|S_1^{(-)}|)$}\,.\end{equation} Moreover from (1.5)
in \cite{ja}, there are finite constants $A$ and $B$ depending only
on the law of $S^{(-)}$ such that for any $a\ge0$
\begin{equation}\label{3}\e(\sigma(a))\le A+Ba\,.\end{equation}
Since $\tilde{\mbox{\em\i}}_0$ is integrable (see \cite{do}), the
result follows from (\ref{1}), (\ref{2}), (\ref{3}) and the
independence between $\tilde{\mbox{\em\i}}_0$ and $S^{(-)}$. \QED

\begin{theorem}\label{main2}
The lower envelope of $X^{(x)}$ at $+\infty$ is described  as follows:\\

\noindent Let $f$ be an increasing function.
\begin{itemize}
\item[$(i)$] If
\[\int^{+\infty}F\left(\frac t{f(t)}\right)\,\frac{dt}t<\infty\,,\]
then for all $\varepsilon>0$, and for all $x\ge0$,
\[\p(X^{(x)}_t<(1-\varepsilon)f(t),\mbox{i.o., as
$t\rightarrow+\infty$})=0\,.\] \item[$(ii)$] If for all $q>0$,
\[\int^{+\infty}F_q\left(\frac t{f(t)}\right)\,\frac{dt}t=\infty\,,\]
then for all $\varepsilon>0$, and for all $x\ge0$,
\[\p(X^{(x)}_t<(1+\varepsilon)f(t),\mbox{i.o., as
$t\rightarrow+\infty$})=1\,.\] \item[$(iii)$] Assume that there
exists $\gamma>1$ such that,
$\limsup_{t\rightarrow+\infty}\p(I>\gamma t)/\p(I>t)<1$. Assume also
that $t\mapsto f(t)/t$ is decreasing. If
\[\int^{+\infty}F\left(\frac t{f(t)}\right)\,\frac{dt}t=\infty\,,\]
then for all $\varepsilon>0$, and for all $x\ge0$,
\[\p(X^{(x)}_t<(1+\varepsilon)f(t),\mbox{i.o., as
$t\rightarrow+\infty$})=1\,.\]
\end{itemize}
\end{theorem}

\noindent{\it Proof}: \noindent We first consider the case where
$x=0$. The proof is very similar to this of Theorem \ref{main1}. We
can follow the proofs of $(i)$, $(ii)$ and $(iii)$ line by line,
replacing the sequences $x_n=r^{-n}$ and $z_n=kr^{-n}$ respectively
by the sequences $x_n=r^n$ and $z_n=kr^n$, and replacing Corollary
\ref{cor1} by Corollary \ref{cor3}. Then with the definition
\[A_n=\{\mbox{There exists $t\in[U(r^n),U(r^{n+1})]$ such that
$X_t^{(0)}<f(t)$.}\}\,,\] we see that the  event $\limsup A_n$
belongs to the tail sigma-field $\cap_t\sigma\{X^{(0)}_s:s\ge t\}$
which is trivial from the representation (\ref{lamp}) and the Markov
property.

The only thing which has to be checked more carefully is the
counterpart at $+\infty$ of the equivalence (\ref{ineq2}). Indeed,
since in that case
$\int_1^\infty\p(rt<f(r^tI)\,dt=\int_{0+}^\infty\p(s/f(s)<I_q,\,s>rI_q)\,ds/(s\log
r)$, in the proof of $(ii)$ and $(iii)$, we need to make sure that
for any $r>1$,
\begin{equation}\label{ineq32}
\int^{+\infty}\p\left(\frac
s{f(s)}<I_q\right)\,\frac{ds}s=+\infty\;\;\;
\mbox{implies}\;\;\;\int^{+\infty}\p\left(\frac
s{f(s)}<I_q<sr\right)\,\frac{ds}s=+\infty\,.\end{equation} To this
aim, note that
\[\int_1^\infty\p\left(\frac s{f(s)}<I_q<sr\right)\,\frac{ds}s=
\int_1^\infty\p\left(\frac s{f(s)}<I_q\right)-\p\left(\frac
s{f(s)}<I_q,\,sr<I_q\right)\,\frac{ds}s\,,\] and since $f$ is
increasing, we have
\[\int_1^\infty\p\left(\frac
s{f(s)}<I_q,\,sr<I_q\right)\,\frac{ds}s<+\infty\;\;\; \mbox{if and
only
if}\;\;\;\int_1^\infty\p\left(s<I_q\right)\,\frac{ds}s<+\infty\,.\]
But
\[\int_1^\infty\p\left(s<I_q\right)\,\frac{ds}s=\e(\log^+I_q)\,.\]
Note that from our hypothesis on $\xi$, we have
$\e(\hat{T}_{-q})<+\infty$, then the conclusion follows from the
inequality
\[\mbox{$\e(\log^+I_q)\le \e\left(\sup_{0\le s\le
\hat{T}_{-q}}\hat{\xi}_s\right)+\e(\hat{T}_{-q})$}\] and Lemma
\ref{finite}. This achieves the proof of the theorem for $x=0$.

Now we prove $(i)$ for any $x>0$. Let $f$ be an increasing function
such that $\int^{+\infty}F\left(\frac
t{f(t)}\right)\,\frac{dt}t<+\infty$. Let $x>0$, put
$S_x=\inf\{t:X^{(0)}_t\ge x\}$ and denote by $\mu_x$ the law of
$X^{(0)}_{S_x}$. From the Markov property at time $S_x$, we have for
all $\varepsilon>0$,
\begin{eqnarray}
&&\p(X^{(0)}_t<(1-\varepsilon)f(t-S_x),\mbox{i.o., as
$t\rightarrow+\infty$})\nonumber\\&=&\int_{[x,\infty)}
\p(X_t^{(y)}<(1-\varepsilon)f(t),\mbox{i.o., as
$t\rightarrow+\infty$})\,\mu_x(
dy)\nonumber\\
&\le&\p(X_t^{(0)}<(1-\varepsilon)f(t),\mbox{i.o., as
$t\rightarrow+\infty$})=0\,.\label{equa1}\end{eqnarray} If $x$ is an
atom of $\mu_x$, then the inequality (\ref{equa1}) shows that
\[\p(X_t^{(x)}<(1-\varepsilon)f(t),\mbox{i.o., as
$t\rightarrow+\infty$})=0\] and the result is proved. Suppose that
$x$ is not an atom of $\mu_x$. Recall from section 1 that
$\log(x_1^{-1}\Gamma)$ is the limit in law of the overshoot process
$\hat{\xi}_{\hat{T}_z}-z$, as $z\rightarrow+\infty$. Moreover, it
follows from \cite{cc}, Theorem 1 that $X^{(0)}_{S_x}\ed \frac
{xx_1}\Gamma$. Hence, from Lemma \ref{overshoot}, we have for any
$\eta>0$, $\mu_x(x,x+\eta)>0$. Then, the inequality  (\ref{equa1})
implies that for any $\eta>0$, there exists $y\in (x,x+\eta)$ such
that $\p(X_t^{(y)}<(1-\varepsilon)f(t),\mbox{i.o., as
$t\rightarrow+\infty$})=0$, for all $\varepsilon>0$. It allows us to
conclude.

Parts $(ii)$ and $(iii)$ can be proved through the same way.\QED

\noindent We recall that to obtain these tests for any scaling index
$\alpha>0$, it suffices to consider the process
$(X^{(0)})^{1/\alpha}$ in the
above theorems. The same remark holds for the results of the next
sections.\\

\section{The regular case}
\setcounter{equation}{0}

\noindent The first type of tail behaviour of $I$ that we consider
is the case where $F$ is regularly varying at infinity, i.e.
\begin{equation}\label{regular}
F(t)\sim \lambda t^{-\gamma}L(t)\,,\;\;\;t\rightarrow+\infty\,,
\end{equation}
where $\gamma>0$ and $L$ is a slowly varying function at $+\infty$.
As shown in the following lemma,  under this assumption, for any
$q>0$ the functions $F_q$ and $F$ are equivalent, i.e. $F_q\asymp
F$.

\begin{lemma}\label{tail} Recall that
$I_q=\int_0^{T_{-q}}\exp\hat{\xi}_s\,ds$ and $F_q(t)=\p(I_q>t)$. If
$(\ref{regular})$ holds then for all $q>0$,
\begin{equation}\label{regular1}
(1-e^{-\gamma q})F(t)\le F_q(t)\le F(t)\,,
\end{equation}
for all $t$ large enough.
\end{lemma}
{\it Proof}: Recall from Lemma \ref{lem3}, that if
$(\hat{\xi}_s,\,s\le \hat{T}_{-q})$ and
$\hat{\xi}'\eqdef(\hat{\xi}_{s+\hat{T}_{-q}}-
\hat{\xi}_{\hat{T}_{-q}},\,s\ge0)$ then
\begin{equation}\label{eqineq}
I=I_q+\exp({\hat{\xi}_{\hat{T}_{-q}}}) I'\le
I_q+e^{-q}I'\end{equation}
 where
$I'=\int_0^\infty\exp\hat{\xi}'_s\,ds$ is a copy of $I$ which is
independent of $I_q$. It yields the second equality of the lemma. To
show the first inequality, we write for all $\delta>0$,
\begin{eqnarray*}
\p(I>(1+\delta)t)&\le&\p(I_q+e^{-q}I'\ge (1+\delta)t)\\
&\le&\p(I_q>t)+\p(e^{-q}I>t)+\p(I_q>\delta t)\p(e^{-q}I>\delta t)\\
&\le&\p(I_q>t)+\p(e^{-q}I>t)+\p(I>\delta t)\p(e^{-q}I>\delta t)\,,
\end{eqnarray*}
so that
\[\liminf_{t\rightarrow+\infty}\frac{\p(I_q>t)}{\p(I>t)}\ge
(1+\delta)^{-\gamma}-e^{-q\gamma}\,,\] and the result follows since
$\delta$ can be chosen arbitrary small.\QED

\noindent The regularity of the behaviour of $F$ allows us to drop
the $\varepsilon$ of Theorems \ref{main1} and \ref{main2} in the
next integral test.

\begin{theorem}\label{regul}
Under condition $(\ref{regular})$, the lower envelope of $X^{(0)}$
at $0$ and at $+\infty$ is as follows:\\

\noindent Let $f$ be an increasing function, such that either
$\lim_{t\downarrow0}f(t)/t=0$, or $\liminf_{t\downarrow0}f(t)/t>0$,
then:
\[\p(X^{(0)}_t<f(t),\mbox{i.o., as
$t\rightarrow0$})=\left\{\begin{array}{l}0\\1\end{array}\right.\,,\]
according as
\[\int_{0+}F\left(\frac t{f(t)}\right)\,\frac{dt}t\;
\left\{\begin{array}{l}<\infty\\=\infty\end{array}\right.\,.\] Let
$g$ be an increasing function, such that either
$\lim_{t\uparrow+\infty}g(t)/t=0$, or
$\liminf_{t\uparrow+\infty}g(t)/t>0$, then for all $x\ge0$,
\[\p(X^{(x)}_t<g(t),\mbox{i.o., as
$t\rightarrow+\infty$})=\left\{\begin{array}{l}0\\1\end{array}\right.\,,\]
according as
\[\int^{+\infty}F\left(\frac t{g(t)}\right)\,\frac{dt}t\;
\left\{\begin{array}{l}<\infty\\=\infty\end{array}\right.\,.\]
\end{theorem}
{\it Proof}: First  let us check that for any constant $\beta>0$:
\begin{equation}\label{equiv}\int_{0+}^\lambda F\left(\frac
s{f(s)}\right)\,\frac{ds}s<\infty\;\;\;\mbox{if and only
if}\;\;\;\int_{0+}^\lambda F\left(\frac {\beta
s}{f(s)}\right)\,\frac{ds}s<\infty\,.\end{equation} From the
hypothesis, either $\lim_{s\downarrow0}f(s)/s=0$, or
$\liminf_{s\downarrow0}f(s)/s>0$. In the first case, we deduce
(\ref{equiv}) from (\ref{regular}). In the second case, since for
any $0<\lambda<\infty$, $\p(I>\lambda)>0$, and
$\limsup_{u\downarrow0}u/{f(u)}<+\infty$, we have for any $s$,
$0<\p\left(\limsup_{u\downarrow0}\frac
u{f(u)}<I\right)<\p\left(\frac s{f(s)}<I\right)$ so both of the
integrals above are infinite.

Now it follows from Theorem \ref{main1} that if
$\int_{0+}F\left(\frac t{f(t)}\right)\,\frac{dt}t<\infty$ then for
all $\varepsilon>0$, $\p(X^{(0)}_t<(1-\varepsilon)f(t),\mbox{i.o.,
as $t\rightarrow0$})=0$. If $\int_{0+}F\left(\frac
t{f(t)}\right)\,\frac{dt}t=\infty$ then from Lemma \ref{tail}, for
all $q>0$, $\int_{0+}F_q\left(\frac
t{f(t)}\right)\,\frac{dt}t=\infty$, and it follows from Theorem
\ref{main1} $(ii)$ that for all $\varepsilon>0$,
$\p(X^{(0)}_t<(1+\varepsilon)f(t),\mbox{i.o., as
$t\rightarrow0$})=1$. Then the equivalence (\ref{equiv}) allows us
to drop $\varepsilon$ in these implications.

The test at $+\infty$ is proven through the same way.\QED

\noindent {\bf Remarks}: \\
{\bf 1.} It is possible to obtain the divergent parts of Theorem
\ref{regul} by applying parts $(iii)$ of Theorems \ref{main1} and
\ref{main2} but then, one has to assume that $f(t)/t$ is an
increasing (respectively a decreasing) function for the test at 0
(respectively at $+\infty$), which is slightly stronger than the
hypothesis on $f$ of Theorem \ref{regul}.\\
{\bf 2.} This result is due to Dvoretzky and Erd\"os \cite{de} and
Motoo \cite{mo} when $X^{(0)}$ is a transient Bessel process, i.e.
the square root of the solution of the SDE:
\begin{equation}\label{Bessel}
Z_t=2\int_0^t\sqrt{Z_s}\,dB_s+\delta t\,,\end{equation} where
$\delta>2$ and $B$ is a standard Brownian motion. (Recall that when
$\delta$ is an integer, $X^{(0)}=\sqrt{Z}$ has the same law as the
norm of the $\delta$-dimensional Brownian motion.) Processes
$X^{(0)}=\sqrt{Z}$ such that $Z$ satisfies the equation
(\ref{Bessel}) with $\delta>2$ are the only continuous self-similar
Markov process with index $\alpha=2$, which drifts towards
$+\infty$. In this particular case, thanks to the time-inversion
property, i.e.:
$$(X_t,\,t>0)\ed(tX_{1/t},\,t>0),$$
we may deduce the test at $+\infty$ from the test at 0.\\
{\bf 3.} A possible way to improve the test at $\infty$ in the
general case (that is in the setting of Theorem \ref{main1}) would
be to first establish it for the Ornstein-Uhlenbeck process
associated to $X^{(0)}$, i.e. $(e^{-t}X^{(0)}(e^t),\,t\ge0)$, as
Motoo did for Bessel processes in \cite{mo}. This would allow us to
consider
test functions which are not necessarily increasing.\\

\noindent {\bf Examples}:\\
{\bf 1.} Examples of such behaviours are given by transient Bessel
processes raised to any power and more generally when the process
$\xi$ satisfies the so called Cramer's condition, that is,
\begin{equation}\label{cramer}
\mbox{there exists $\gamma>0$ such that $E(e^{-\gamma\xi_1})=1$.}
\end{equation}
In that case, Rivero \cite{ri2} and Maulik and Zwart \cite{mz}
proved by using results of Kesten and Goldie on tails of solutions
of random equations that the behavior of $P(I>t)$ is given by
\begin{equation}\label{cramer1}
F(t)\sim Ct^{-\gamma}\,,\;\;\;\mbox{as $t\rightarrow+\infty$,}
\end{equation}
where the constant $C$ is explicitly computed in \cite{ri} and
\cite{mz}.\\

\noindent {\bf 2.} Stable L\'evy processes conditioned to stay
positive are themselves positive self-similar Markov processes which
belong to the regular case. These processes are defined as
$h$-processes of the initial process when it starts from $x>0$ and
killed at its first exit time of $(0,\infty)$. Denote by $(q_{t})$
the semigroup of a stable L\'evy process $Y$ with index
$\alpha\in(0,2]$, killed at time $R=\inf\{t:Y_t\le0\}$. The function
$h(x)=x^{\alpha(1-\rho)}$, where $\rho=\p(Y_1\ge0)$, is invariant
for the semi-group $(q_{t})$, i.e. for all $x\ge0 $ and $t\ge 0$,
$E_x(h(Y_{t})\mbox{\rm 1\hspace{-0.033in}I}_{\{t<R\}})=h(x)$, ($E_x$
denotes the law of $Y+x$). The L\'{e}vy process $Y$ conditioned to
stay positive is the strong Markov process whose semigroup is
\begin{equation}  \label{pf}
p_{t}^{\uparrow }(x,dy):=\frac{h(y)}{h(x)}q_{t}(x,dy),\mbox{
}x>0,y>0,\mbox{ }t\geq 0\,.
\end{equation}
We will denote this process by $X^{(x)}$ when it is issued from
$x>0$. We refer to \cite{ch} for more on the definition of L\'evy
processes conditioned to stay positive and for a proof of the above
facts. It is  easy to check that the process $X^{(x)}$ is
self-similar and drifts towards $+\infty$. Moreover, it is proved in
\cite{ch}, Theorem 6 that $X^{(x)}$ converges weakly as
$x\rightarrow0$ towards a non degenerated process $X^{(0)}$ in the
Skorohod's space, so from \cite{cc}, the underlying L\'evy process
in the Lamperti representation of $X^{(x)}$ satisfies condition
$(H)$.

We can check that the law of $X^{(x)}$ belongs to the regular case
by using the equality in law (\ref{lamp2}). Indeed, it follows from
Proposition \ref{ret} and Theorem 4 in \cite{ch} that the law of the
exponential functional $I$ is given by
\begin{equation}\label{loiI}
\p(t<x^\alpha
I)=x^{1-\alpha\rho}E_{-x}(\hat{Y}_t^{\alpha\rho-1}\ind_{\{t<\hat{R}\}})\,,
\end{equation}
where $\hat{Y}=-Y$ and $\hat{R}=\inf\{t:\hat{Y}_t\le0\}$. If $Y$
(and thus $X^{(0)}$) has no positive jumps, then $\alpha\rho=1$ and
it follows from (\ref{loiI}) and Lemma 1 in \cite{ch2} that
\begin{equation}\label{tailI}
\p(t<I)=Ct^{-\rho}\,.
\end{equation}
We conjecture that (\ref{tailI}) is also valid when $Y$ has positive
jumps. We also emphasize the possibility that the underlying L\'evy
process in the Lamperti representation of $X^{(x)}$ even satisfies
(\ref{cramer}) with $\gamma=\rho$.\\

\section{The log regular case}
\setcounter{equation}{0}

The second type of behaviour that we shall consider is when $\log F$
is regularly varying at $+\infty$, i.e.
\begin{equation}\label{reg}
-\log F(t)\sim \lambda t^{\beta}L(t)\,,\;\;\;\mbox{as
$t\rightarrow\infty$,}
\end{equation}
where $\lambda>0$, $\beta>0$ and  $L$ is a function which varies
slowly at $+\infty$. Define the function $\psi$ by
\begin{equation}\label{def}\psi(t)\eqdef \frac{t}{\inf\{s:1/F(s)>|\log
t|\}}\,,\;\;\;t>0\,,\;t\neq1\,.\end{equation} Then the lower
envelope of $X^{(0)}$ may be described as follows:

\begin{theorem}\label{logregul}
Under condition $(\ref{reg})$, the process $X^{(0)}$ satisfies the
following law of the iterated logarithm:
\begin{itemize}\item[$(i)$]
\begin{equation}
\liminf_{t\rightarrow0}\frac{X^{(0)}_t}{\psi(t)}=1\,,\;\;\mbox{almost
surely.}
\end{equation}
\item[$(ii)$] For all $x\ge0$,
\begin{equation}
\liminf_{t\rightarrow+\infty}\frac{X^{(x)}_t}{\psi(t)}=1\,,\;\;\mbox{almost
surely.}
\end{equation}
\end{itemize}
\end{theorem}
{\it Proof}: We shall apply Theorem \ref{main1}. We first have to
check that under hypothesis (\ref{reg}), the conditions of part
$(iii)$ in Theorem  \ref{main1} are satisfied. On the one hand, from
(\ref{reg}) we deduce that for any $\gamma>1$, $\limsup F(\gamma
t)/F(t)=0$. On the other hand, it is easy to see that both $\psi(t)$
and $\psi(t)/t$ are increasing in a neighbourhood of 0.

Let $\overline{L}$ be a slowly varying function such that
\begin{equation}\label{inv}
-\log F(\lambda^{-1/\beta} t^{1/\beta} \overline{L}(t))\sim
t\,,\;\;\;\mbox{as $t\rightarrow+\infty$.}
\end{equation}
Th. 1.5.12, p.28 in \cite{bgt} ensures that such a function exists
and that
\begin{equation}\label{inv2}
\inf\{s:-\log F(s)>t\}\sim\lambda^{-1/\beta} t^{1/\beta}
\overline{L}(t)\,,\;\;\;\mbox{as $t\rightarrow+\infty$.}
\end{equation}
Then we have for all $k_1<1$ and $k_2>1$ and for all $t$
sufficiently large,
\begin{eqnarray*}
k_1\lambda^{-1/\beta} t^{1/\beta} \overline{L}(t)\le \inf\{s:-\log
F(s)>t\}\le k_2\lambda^{-1/\beta} t^{1/\beta} \overline{L}(t)
\end{eqnarray*}
so that for $\psi$ defined above and for all $k_2'>0$,
\begin{equation}\label{ine1}
-\log F\left(\frac{t\,k'_2}{k_2\psi(t)}\right)\le -\log
F(k_2'\lambda^{-1/\beta}(\log|\log t|)^{1/\beta}
\overline{L}(\log|\log t|))
\end{equation}
for all $t$ sufficiently small. But from (\ref{inv}), for all
$k''_2>1$ and for all $t$ sufficiently small,
\[-\log
F(k_2'\lambda^{-1/\beta}(\log|\log t|)^{1/\beta}
\overline{L}(\log|\log t|))\le k''_2{k'_2}^\beta\log|\log t|\,,\]
hence
\[F\left(\frac{t\,k'_2}{k_2\psi(t)}\right)\ge  (|\log
t|)^{-k''_2{k'_2}^\beta}\,.\] By choosing $k'_2<1$ and
$k_2''<(k'_2)^{-\beta}$, we obtain the convergence of the integral
\[\int_{0+}F\left(\frac{t\,k'_2}{k_2\psi(t)}\right)\,\frac{dt}t\,,\]
for all $k_2>1$ and $k_2'<1$, which proves that for all
$\varepsilon>0$,
\[\p(X^{(0)}_t<(1+\varepsilon)\psi(t),\mbox{i.o., as
$t\rightarrow0$})=1\] from Theorem \ref{main1} $(iii)$. The
convergent part is proven through the same way so that from Threorem
\ref{main1} $(i)$, one has for all $\varepsilon>0$,
\[\p(X^{(0)}_t<(1-\varepsilon)\psi(t),\mbox{i.o., as
$t\rightarrow0$})=0\] and the conclusion follows.

Condition (\ref{reg}) implies that $\psi(t)$ is increasing in a
neighbourhood of $+\infty$ whereas $\psi(t)/t$ is decreasing in a
neighbourhood of $+\infty$. Hence, the proof of the result at
$+\infty$ is done through the same way as at 0, by using Theorem
\ref{main2}, $(i)$ and $(iii)$. \QED

\noindent {\bf Example}:\\
 An example of such a behaviour is provided
by the case where the process $X^{(0)}$ is increasing, that is when
the underlying L\'evy process $\xi$ is a subordinator. Then Rivero
\cite{ri}, see also Maulik and Zwart \cite{mz} proved that when the
Laplace exponent $\phi$ of $\xi$ which is defined by
\[\exp(-t\phi(\lambda))=\e(\exp(\lambda\hat{\xi}_t))\,,\;\;\;\lambda>0,\;\;t\ge0\]
is regularly varying at $+\infty$ with index $\beta\in(0,1)$, the
upper tail of the law
 of $I$ and the asymptotic behavior of $\phi$ at $+\infty$ are related as
 follows:
\begin{proposition}\label{propri}
Suppose that $\xi$ is a subordinator whose Laplace exponent $\phi$
varies regularly at infinity with index $\beta\in(0,1)$, then
\[-\log F(t)\sim (1-\beta)\phi^{\leftarrow}(t)\,,\;\;\;\mbox{as
$t\rightarrow\infty$,}\] where
$\phi^{\leftarrow}(t)=\inf\{s>0:s/\phi(s)>t\}$.
\end{proposition}
Then by using an argument based on the study of the associated
Ornstein-Uhlenbeck process $(e^{-t}X^{(0)}(e^t),\,t\ge0)$ Rivero
\cite{ri} derived from Proposition \ref{propri} the following
result. Define
\[\varphi(t)=\frac{\phi(\log|\log t|)}{\log|\log
t|},\;\;\;t>e\,.\]
\begin{corollary}\label{corri}
If $\phi$ is regularly varying at infinity with index
$\beta\in(0,1)$ then
\[\liminf_{t\downarrow0}\frac{X^{(0)}}{t\varphi(t)}=(1-\beta)^{1-\beta}\;\;\;\mbox{and}
\;\;\;\liminf_{t\uparrow+\infty}\frac{X^{(0)}}{t\varphi(t)}=(1-\beta)^{1-\beta},\;\;\;a.s.\]
\end{corollary}
This corollary is also a consequence of Proposition \ref{propri} and
Theorem \ref{logregul}. To establish Corollary \ref{corri}, Rivero
assumed moreover that the density of the law of the exponential
functional $I$ is decreasing and bounded in a neighbourhood of
$+\infty$. This additional assumption is actually needed to
establish an integral test which involves the density of
$I$ and which implies Corollary \ref{corri}.\\

\vspace*{0.5in}

\noindent \textbf{Acknowledgment}: We are much indebted to professor
Zhan Shi for many fruitful discussions on integral tests and laws of
the iterated logarithm. We also warmly thank professor Ron Doney for
having provided to us the idea of the proof of Lemma \ref{finite}.

\end{document}